\documentclass[12pt]{amsart}
\usepackage{enumerate}
\usepackage{graphicx,graphics}
\usepackage{soul}
\usepackage{amsfonts}
\usepackage{amssymb}
\usepackage{amsthm}
\usepackage{amsmath}
\usepackage{marginnote}
\usepackage{mathrsfs,color}
\usepackage{tikz} % diagramas
\usepackage{hyperref}
\input{xy}
\xyoption{all}

\numberwithin{equation}{section}

\newtheorem{theorem}{Theorem}[section]
\newtheorem{corollary}{Corollary}[theorem]
\newtheorem{lemma}[theorem]{Lemma}

\theoremstyle{definition}
\newtheorem{definition}[theorem]{Definition}
\theoremstyle{definition}
\newtheorem{problem}[theorem]{Problem}
\theoremstyle{definition}
\newtheorem{example}[theorem]{Example}
\theoremstyle{definition}

\theoremstyle{definition}
\newtheorem{question}[theorem]{Question}
\newcommand*{\QEDA}{\hfill\ensuremath{\blacksquare}}
\begin{document}

\title[Equivariant means]{Equivariant means}

\author{Natalia Jonard-P\'erez}\email{nat@ciencias.unam.mx}

%\author[N. Jonard-P\'erez, A. L\'opez-Poo]{Natalia Jonard-P\'erez and Ananda L\'opez-Poo}

\author{Ananda L\'opez-Poo}
\email{anandalc95@gmail.com}

\address{Departamento de  Matem\'aticas,
Facultad de Ciencias, Universidad Nacional Aut\'onoma de M\'exico, 04510 Ciudad de M\'exico, M\'exico.
            }


  % Social Choice
%01A60, %History of mathematics in the 20th century 
%55P20, %Eilenberg-Mac Lane spaces

\thanks{This work has been supported by  PAPIIT grant IN101622  (UNAM, M\'exico). Also, the second author has been supported by Conahcyt grant 712523.}

\begin{abstract} 
An $n$-mean (also called  a ``topological social choice rule'') on a topological space $X$ is a continuous function $p:X^n\to X$ satisfying $p(x,\dots, x)=x$ for every $x\in X$ and $p(x_1,\dots, x_n)=p(x_{\sigma(1)},\dots x_{\sigma(n)})$  for any permutation $\sigma$ of $\{1,\dots, n\}$. If, in addition, $X$ is a $G$-space and $p$ is equivariant with respect to the diagonal action of $G$ on $X^n$, we say that $p$ is an equivariant $n$-mean. 
In this paper, we continue the work initiated by H. Ju\'arez-Anguiano about conditions on a $G$-space $X$, under which the existence of an equivariant $n$-mean guarantees that $X$ is a $G$-AR. We also explore this problem when we remove the symmetry condition on the definition of an $n$-mean.  
\end{abstract}

\maketitle

\section{Introduction}

Let $X$ be a topological space. Given a positive integer $n$, an $n$-\textit{mean} in $X$ is a continuous function $p:X^n\to X$  satisfying the following two conditions:

\begin{enumerate} 
\item[(M1)] $p\left(x,\ldots, x\right)=x$ for every $x\in X$,
\item[(M2)] $p\left(x_{1},\ldots,x_{n}\right)=p\left(x_{\sigma\left(1\right)},\ldots,x_{\sigma\left(n\right)}\right)$ for each $\left(x_{1},\ldots,x_{n}\right)\in X^{n}$ and every permutation $\sigma$ of $\left\{1,\ldots,n\right\}$.
\end{enumerate} 

The arithmetic mean of $n$ real numbers is a simple example of an $n$-mean in $\mathbb R$. In the same way,  the harmonic mean in $(0,\infty)$ or the geometric mean in $[0,\infty)$ are classical examples of such functions. Notice that in all these examples the domain is always a convex set (and therefore it is contractible).

In 1954,  B. Eckmann proved that if $X$ is a compact connected polyhedron, then the existence of an $n$-mean in $X$ (for every $n\in\mathbb N$) implies that $X$ is contractible (\cite{eck}). Some years later,  B. Eckmann, T. Ganea and P. J. Hilton improved that result by showing that if $X$ is a connected topological space with the homotopy type of a polyhedron with finitely generated homology groups and vanishing homology above a certain dimension, then the existence of an $n$-mean in $X$ (for some $n\geq 2$) is equivalent to $X$ being contractible (\cite[Corollary 3.2 and p. 90]{eck2}). As  was observed by H. Ju\'arez-Anguiano in \cite{hugo}, since every ANR is homotopically equivalent to a polyhedron and since the class of contractible ANR-spaces is precisely the class of AR-spaces, then the result obtained by B. Eckmann, T. Ganea and P. J. Hilton can be stated as follows.  

 \begin{theorem} \label{t:caracterizacion n mean original}
Let $X$ be a connected ANR with finitely generated homology groups such that almost all vanish. Then the following statements are equivalent.

\begin{enumerate}
\item[(1)] There exists an $n$-mean $p:X^{n}\rightarrow X$ for every $n\geq 2$.
\item[(2)] There exists an $n$-mean $p:X^{n}\rightarrow X$ for some $n\geq 2$.
\item[(3)] $X$ is an AR.
\end{enumerate} 
\end{theorem}

Inspired by Arrow's imposibility paradox (see, e.g. \cite{Arrow}), economists G. Chichilnisky and G. Heal rediscovered in 1983 implications $(1)\Longleftrightarrow (3)$   (without knowing Eckmann's result), but under a different terminology and motivation (\cite{ch}).
In their work, the topological space $X$ was a connected $CW$-complex called a space of \textit{preferences}, conditions (M1) and (M2) were called \textit{unanimity} and \textit{anonymity}, respectively, and an $n$-mean $p:X^n\to X$ was called a \textit{topological social choice rule}. 

Twenty years later, following the work and terminology of G.Chichilnisky and G. Heal (and also without knowing Eckmann's result), S. Weinberger rediscovered implication $(2)\Longrightarrow (3)$  (\cite[Theorem 1.1]{hip}), among other results.
In his paper, Weinberger shows that if a non-contractible connected $CW$-complex $X$ admits an $n$-mean $p:X^n\to X$ for every $n\in\mathbb N$, then these means will misbehave. More precisely, for a suitable metric $d$ in $X$ and for every $K>0$, there exist points $x_1,\dots, x_n\in X$ such that $d\left(p(x_{1},\ldots,x_{n}),x_{i}\right)>K$ for all $i=1,\ldots,n$. Weinberger called this kind of function \textit{solomonic} (see \cite[Proposition 3.2]{hip} and \cite[\S 5.4]{eck3}). 
%We refer the reader to the little survey \cite{eck3} for more insight into the theory of $n$-means.

Another twenty years later, H. Juárez-Anguiano proved an equivariant version of Theorem~\ref{t:caracterizacion n mean original} for the case when $X$ is a $G$-space, $G$ is a compact group and the $n$-means are $G$-equivariant\footnote{We refer the reader to Section 2 for the definition of an equivariant mean and any other technical notion.}.  Later, in \cite[Theorem 4.5]{hugo2}, Juárez-Anguiano improved this result for the case when the group $G$ is almost-connected and $X$ is a proper $G$-space (see Theorem~\ref{hug}). In both theorems, an extra hypothesis was needed: that for every closed subgroup $H\leq G$, the set of all $H$-fixed points $X^H$ is connected and has finitely generated homology groups such that almost all vanish.  For this reason, in \cite{hugo}, H. Juárez-Anguiano also posed the following question.

\begin{question} \label{preguntahugo}
Let $G$ be a compact group and let $X$ be a compact and connected metrizable $G$-space that is a $G$-ANR. If there exists an equivariant $n$-mean $p:X^{n}\rightarrow X$ for some $n\geq 2$, then is $X$ a $G$-AR?
\end{question} 

This question adds to the list of open problems that seek to characterize when a $G$-space is a $G$-AR. One of those open problems is the following.

\begin{problem}[J. Jaworowski's problem] \label{jaw}
Let $G$ be a compact Lie group and $X$ be a metrizable $G$-space that has a finite number of orbit types (for example, when $G$ is a finite group). Assume that for every closed subgroup $H$ of $G$ the set $X^{H}$ is an AR. Is $X$ a $G$-AR?
\end{problem}
Despite many efforts to solve it (see \cite{s4, we, ww, wo}), this problem remains open even in the simplest case, when $X$ is homeomorphic to the Hilbert cube $Q=[0,1]^{\mathbb N}$, $G$ is the group $\mathbb{Z}_{2}$ and $X^{\mathbb{Z}_{2}}$ is a singleton\footnote{This particular case of Jaworowski's problem is equivalent to another open problem known as Anderson's problem (see \cite[Subsection 3.2]{s2}).}.

\medskip

Inspired by problems~\ref{preguntahugo} and \ref{jaw}, in this paper we continue the work initiated in \cite{hugo} by investigating when does the existence of an equivariant $n$-mean guarantee that the space is a $G$-AR. 

The paper is organized as follows. In Section 2, we introduce all the basic notions and results regarding $G$-spaces and equivariant $n$-means. In Section 3, we explore the case when the group $G$ is finite. We prove that under the hypotheses of Jaworowski's problem, the existence of a suitable equivariant $n$-mean guarantees that the space is a $G$-AR (Theorem~\ref{t:answer to Jaworowski}). We also give a partial answer to Question~\ref{preguntahugo} for the case when $G$ is finite
(Theorem~\ref{t:n mean implica G-ANR} and Corollary~\ref{c:problema hugo}).

After Weinberger's result about solomonic means, in Section 4 we investigate Question~\ref{preguntahugo} (and Theorem~\ref{t:caracterizacion n mean original}) from a different perspective. In order to explain it, first let us point out that from the point of view of Social Choice Theory, the axiom of anonymity (M2) in the definition of an $n$-mean is a way to avoid that the social choice rule (i.e. the $n$-mean) becomes dictatorial. A map $p:X^n\to X$ is \textit{dictatorial} if it is the projection on a fixed coordinate $j\in\{1,\dots, n\}$ (in this case, the coordinate $j$ represents the idea of a dictator whose will is the only possible output of the map $p$). However, there are other ways to ensure that a map $p:X^n\to X$ is not dictatorial. 

For instance, since means in a non-contractible space are solomonic,  if we want that the existence of an $n$-mean guarantees the contractibility of the space, the first thing to ask is that the $n$-mean is not solomonic. One way to do this is by considering a restriction on a given distance in the space. 
Hence, in Definition~\ref{d:quasimean} we introduce the notion of a contractive $n$-quasi-mean on a metric space $(X,d)$. We show that those maps satisfy property (M1). And even if they do not satisfy axiom (M2), they are not dictatorial. 
Then we prove that in a complete metric space, the existence of a contractive $n$-quasi-mean implies that the space is contractible. And, if additionally $X$ is a $G$-ANR, the $n$-quasi-mean is equivariant, and the set of $G$-fixed points is not empty, then $X$ is a $G$-AR (Corollary~\ref{contc}).

\section{Preliminaries}

Throughout this paper, the letter $G$ will always denote a compact Hausdorff topological group, and we will use the word \text{map} as a synonym of a continuous function.

An \textit{action} of a group $G$ on a topological space $X$ is a continuous function $\alpha:G\times X\rightarrow X$ such that $\alpha\left(e,x\right)=x$ and $\alpha\left(g,\alpha\left(h,x\right)\right)=\alpha\left(gh,x\right)$, where $g,h\in G$, $x\in X$ and $e$ is the identity element of $G$. To simplify the notation, we write $gx$ instead of $\alpha\left(g,x\right)$. A topological space $X$ equipped with an action of a topological group $G$ is called a \textit{$G$-space}. We refer the reader to \cite{Bredon} for basic notions of the theory of G-spaces. However, we recall here some special definitions and results that will be used throughout this text.

Let $X$ be a $G$-space. If $A$ is a subset of $X$, the \textit{saturation} of $A$, denoted by $G\left(A\right)$, is defined as the set $\left\{ga \mid g\in G, \mbox{ } a\in A\right\}$. If $G\left(A\right)=A$ we say that $A$ is \textit{invariant}. If $A=\left\{x\right\}$ for some $x\in X$, we simply write $G\left(x\right)$ instead of $G\left(\left\{x\right\}\right)$. In this case, we say that $G\left(x\right)$ is the \textit{orbit} of $x$.

If $x\in G$, we denote by $G_{x}$ the \textit{stabilizer} of $x$, defined by $G_{x}:=\left\{g\in G \mid gx=x\right\}$. Given a subgroup $H$ of $G$, we denote the set of $H$-fixed points of $X$ by $X^{H}$. Namely,  
\[
X^{H}:=\left\{x\in X \mid hx=x \mbox{ for each } h\in H\right\}.
\]
If $x\in X^{G}$, we say that $x$ is a \textit{$G$-fixed point}. It is well known that, for every subgroup $H\leq G$, the set $X^H$ is always closed. 

The family of all subgroups of $G$ that are conjugated with $H$ will be denoted by $(H)$. That is, $(H)=\left\{gHg^{-1} \mid g\in G\right\}$. For every $x\in X$, $(G_{x})$ will be called the \textit{orbit type} of the orbit $G\left(x\right)$. If $(H)$ is the orbit type of some orbit of $X$, then $(H)$ will be called an \textit{orbit type} of $X$.

A function $f:X\rightarrow Y$ between $G$-spaces is called $G$-\textit{equivariant} (or simply \textit{equivariant}) if $f\left(gx\right)=gf\left(x\right)$ for every $g\in G$, $x\in X$. If in addition $f$ is continuous, then we say that $f$ is a $G$-\textit{map}.

A metrizable $G$-space $Y$ is called a \textit{$G$-absolute neighborhood extensor} (denoted by $G$-ANE) if for every closed and invariant subset $A$ of a metrizable $G$-space $X$ and every $G$-map $f:A\rightarrow Y$, there exist an invariant neighborhood $U$ of $A$ in $X$ and a $G$-map $F:U\rightarrow Y$ such that $F\restriction_{A}=f$. If we can always take $U=X$, then we say that $Y$ is a \textit{$G$-absolute extensor} (denoted by $G$-AE).

Similarly, a metrizable $G$-space $Y$ is a \textit{$G$-absolute neighborhood retract} (denoted by $G$-ANR) if for every metrizable $G$-space $X$ that contains $Y$ as a closed and invariant subset there exist an invariant neighborhood $U$ of $Y$ in $X$ and an equivariant retraction $r:U\rightarrow Y$. If we can always take $U=X$,  we say that $Y$ is a \textit{$G$-absolute retract} (denoted by $G$-AR).

If we take $G=\left\{e\right\}$ in the previous definitions, we obtain the classical notions of \textit{absolute neighborhood extensor} (ANE), \textit{absolute extensor} (AE), \textit{absolute neighborhood retract} (ANR) and \textit{absolute retract} (AR).

In \cite[Theorem 14]{s1} it was proved that a metrizable $G$-space $X$ is a $G$-ANE ($G$-AE) if and only if it is a $G$-ANR ($G$-AR).

Let $f,g:X\rightarrow Y$ be $G$-maps between $G$-spaces. A \textit{$G$-homotopy} between $f$ and $g$ is a homotopy $H:X\times \left[0,1\right]\rightarrow Y$ such that $H_{0}=f$, $H_{1}=g$ and $H_{t}$ is a $G$-map for each $t\in \left[0,1\right]$. A $G$-space $X$ is called \textit{$G$-contractible} if there is a $G$-homotopy between the identity map of $X$ and a constant function $X\rightarrow \left\{x_{0}\right\}$ (notice that in this case $x_{0}$ must be  a $G$-fixed point). An invariant subset $A$ of a $G$-space $X$ is a \textit{$G$-strong deformation retract} of $X$ if there exists a $G$-homotopy $H:X\times \left[0,1\right]\rightarrow Y$ such that $H_{0}$ is the identity map of $X$, $H_{1}$ is a retraction from $X$ onto $A$ and $H\left(a,t\right)=a$ for all $a\in A$ and $t\in \left[0,1\right]$.

In the next section we will use the following theorems due to S. Antonyan.
\begin{theorem} \cite[Theorem 6]{s5} \label{gc}
Let $X$ be a metrizable G-space. Then $X$ is a $G$-AR if and only if $X$ is a $G$-ANR and is $G$-contractible.
\end{theorem}

\begin{theorem} \cite[Theorem 3.7]{s3} \label{t:Sergey strong deformation retract}
Let $G$ be a compact Lie group. A metrizable $G$-space $X$ is a $G$-AR if and only if $X$ is an AR and, for every closed subgroup $H$ of $G$, $X^{H}$ is a $H$-strong deformation retract of $X$.
\end{theorem}

An \textit{$n$-mean} on a topological space $X$, is a continuous function $p:X^{n}\rightarrow X$ that satisfies the following conditions:
\begin{enumerate} 
\item[(M1)] $p\left(x,\ldots, x\right)=x$ for every $x\in X$,
\item[(M2)] $p\left(x_{1},\ldots,x_{n}\right)=p\left(x_{\sigma\left(1\right)},\ldots,x_{\sigma\left(n\right)}\right)$ for each $\left(x_{1},\ldots,x_{n}\right)\in X^{n}$ and every permutation $\sigma$ of $\left\{1,\ldots,n\right\}$.
\end{enumerate} 
When $n=2$, we simply say that $p$ is a \textit{mean}.

\begin{definition}
    Let $X$ be a topological space. An  \textit{$n$-quasi-mean} is a continuous function $p:X^{n}\rightarrow X$ satisfying  condition (M1). If, in addition, $n=2$, we simply say that $p$ is a \textit{quasi-mean}.
\end{definition}

An \textit{equivariant $n$-mean} (\textit{equivariant $n$-quasi-mean}) on a $G$-space $X$ is an $n$-mean ($n$-quasi-mean) $p:X^{n}\rightarrow X$ that additionally satisfies $p\left(gx_{1},\ldots,gx_{n}\right)=gp\left(x_{1},\ldots,x_{n}\right)$ for every $\left(x_{1},\ldots,x_{n}\right)\in X^{n}$ and $g\in G$. 
When $n=2$, we simply say that $p$ is an \textit{equivariant mean} (\textit{equivariant quasi-mean}).

The following is the most general equivariant version of Theorem~\ref{t:caracterizacion n mean original}.\footnote{In \cite[Theorem 3.3]{hugo} the reader can find a less general version for the case where the group $G$ is compact.}

\begin{theorem} (\cite[Theorem 4.5]{hugo2})\label{hug}
Let $G$ be an almost-connected group and $X$ be a proper $G$-space. Assume that $X$ is a $G$-ANR and  that for each compact subgroup $H$ of $G$, $X^{H}$ is connected and has finitely generated homology groups such that almost all vanish. Then the following conditions are equivalent.

\begin{enumerate} 
\item[1)] There exists an equivariant $n$-mean $p:X^{n}\rightarrow X$ for every $n\geq 2$.
\item[2)] There exists an equivariant $n$-mean $p:X^{n}\rightarrow X$ for some $n\geq 2$.
\item[3)] $X$ is a $G$-AR.
\end{enumerate} 
\end{theorem}

We omit the definitions of an almost connected group and a proper $G$-space $X$, since we do not need them. However, the interested reader can find them in \cite[p. 3]{hugo2}. We only point out that if the group $G$ is compact, then $G$ is almost connected and every $G$-space is proper.

\section{Equivariant means for finite groups}

In this section, we provide a positive answer for Jaworowski's problem in the case when $G$ is finite and there exists an equivariant $|G|$-mean. 
Under these two hypotheses, we also give a partial answer to Question~\ref{preguntahugo}. 

Let us start by noticing that the existence of an $n$-mean ($n$-equivariant mean) guarantees the existence of $k$-means  ($k$-equivariant means) for every divisor $k$ of $n$, as we prove in the following lemma. 

\begin{lemma} \label{l: |G| divides n}
 If $p:X^{n}\rightarrow X$ is an $n$-mean and $k\in \mathbb{N}$ divides $n$,   then there exists a $k$-mean $q:X^{k}\rightarrow X$. If in addition $p$ is an equivariant $n$-mean, then $q$ can be chosen to be equivariant too. 
\end{lemma}

\begin{proof}
Define  $q:X^{k}\rightarrow X$ by
 \[
 q\left(x_{1},\ldots,x_{k}\right)= p\left(x_{1},\ldots,x_{k},\ldots,x_{1},\ldots,x_{k}\right),
 \] 
 where $x_{1},\ldots,x_{k}$ are repeated $\frac{n}{k}$ times in the argument of $p$. Clearly $q$ is a $k$-mean. And if $p$ is equivariant, then $q$ is equivariant too.
\end{proof}

\begin{lemma}\label{l:n mean +homotopy=G-homotopy} Let $G$ be a finite group and 
    let $X$ be a $G$-space.  Assume that there exists an equivariant $n$-mean $p:X^n\to X$, where $n=|G|$. If  $\Phi:X\times[0,1]\to X$ is a homotopy and $g_1,\dots, g_n$ are the elements of $G$, then the map $\Psi:X\times \left[0,1\right]\rightarrow X$ defined by
\[
\Psi\left(x,t\right):=p\left(g_{1}^{-1}\Phi\left(g_{1}x,t\right),\ldots,g_{n}^{-1}\Phi\left(g_{n}x,t\right)\right),
\]
is a $G$-homotopy. Furthermore, 
\begin{enumerate}[\rm(1)]
    \item If $\Phi_0$ is the identity map on $X$, then $\Psi_0$ is the identity map on $X$.
    \item If $\Phi_1$ is a constant map, then $\Psi_1$ is a constant map. 
    \item If $\Phi$ is a contraction, then $\Psi$ is an equivariant contraction. 
\end{enumerate}  
    \end{lemma}  

\begin{proof}
Let us start by proving that $\Psi$ is a $G$-homotopy. Clearly, $\Psi$ is continuous. 
Observe that for every $h\in G$, the rule $g_i\to g_{i}h$ defines a bijection on $G$. Then, since $p$ is equivariant and satisfies property (M2), we have that 
\begin{align*}
\Psi(hx,t)&=p\big(g_1^{-1}\Phi(g_{1}hx,t),\dots, g_{n}^{-1}\Phi(g_nhx,t)\big)\\
&=hh^{-1}p\big(g_1^{-1}\Phi(g_{1}hx,t),\dots, g_{n}^{-1}\Phi(g_nhx,t)\big)\\
&=hp\big(h^{-1}g_1^{-1}\Phi(g_{1}hx,t),\dots, h^{-1}g_{n}^{-1}\Phi(g_nhx,t)\big)\\
&=hp\big((g_1h)^{-1}\Phi(g_{1}hx,t),\dots, (g_{n}h)^{-1}\Phi(g_nhx,t)\big)\\
&=hp\big(g_1^{-1}\Phi(g_{1}x,t),\dots, g_{n}^{-1}\Phi(g_nh,t)\big)\\
&=h\Psi(x,t)
\end{align*}
for every $h\in G$ and $(x,t)\in X\times [0,1]$. Hence, $\Psi$ is a $G$-homotopy.

(1) Assume that $\Phi(x,0)=x$ for all $x\in X$. Then, by condition (M1) we have that
\begin{align*}
\Psi\left(x,0\right)&=p\left(g_{1}^{-1}\Phi\left(g_{1}x,0\right),\ldots,g_{n}^{-1}\Phi\left(g_{n}x,0\right)\right)\\
&=p\left(g_{1}^{-1}g_{1}x,\ldots,g_{n}^{-1}g_{n}x\right)\\
&=x.
\end{align*}

(2) If $c\in X$ is such that $\Phi(x, 1)=c$ for every $x\in X$, then 
\begin{align*}
\Psi\left(x,1\right)&=p\left(g_{1}^{-1}\Phi\left(g_{1}x,1\right),\ldots,g_{n}^{-1}\Phi\left(g_{n}x,1\right)\right)\\
&=p\left(g_{1}^{-1}c,\ldots,g_{n}^{-1}c\right).
\end{align*}
Hence, $\Psi_1$ is a constant map.

(3) Follows directly from (1) and (2).
\end{proof}
    
\begin{lemma}\label{l:fixed points are strong deformation retracts}

Let $X$ be a metrizable $G$-space, where $G$ is a finite group, and let $H\leq G$ be an arbitrary subgroup with $|H|=n$. Assume that there exists an equivariant $n$-mean $p:X^{n}\rightarrow X$. Then the following hold:
\begin{enumerate}[\rm(1)]
    \item The set $X^{H}$ is nonempty.
    \item If $X$ is an AR and $X^{H}$ is a retract of $X$, then $X^H$ is  an $H$-strong deformation retract of $X$.
    \end{enumerate}
\end{lemma}

\begin{proof}
Since $|H|=n$, we can assume that $H=\{g_1,\dots, g_n\}$.

(1) Let $x\in X $ be an arbitrary point and define
$$x_0:=p(g_1x, \dots, g_nx).$$
We claim that $x_0\in X^H$. Indeed, if $h\in H$, observe that $g_i\to hg_i$ defines a bijection on $H$. Then, using the fact that
 $p$ is equivariant and satisfies condition (M2), we conclude that
\begin{align*}
    hx_0=&hp(g_1x,\dots, g_nx)=p(hg_1x, \dots, hg_nx)\\
    &=p(g_1x,\dots, g_nx)=x_0.
\end{align*}
 This implies that $x_0\in X^H$, as desired. 
 
(2) Let $r:X\rightarrow X^{H}$ be a retraction. Consider $A=\left(X\times \left\{0\right\}\right)\cup \left(X\times \left\{1\right\}\right)\cup \left(X^{H}\times \left[0,1\right]\right)$ and define the map $\varphi:A\rightarrow X$ given by
\[
\varphi\left(x,t\right):=
\left\{\begin{array}{ll} 
x & \mbox{ if } \left(x,t\right)\in X\times \left\{0\right\},\\
r\left(x\right) & \mbox{ if } \left(x,t\right)\in \left(X\times \left\{1\right\}\right)\cup \left(X^{H}\times \left[0,1\right]\right).\end{array}
\right. 
\]
Since $X$ is a metrizable AR, it is also an AE. Then, since $A$ is closed in $X\times \left[0,1\right]$, there exists an extension $\Phi:X\times \left[0,1\right]\rightarrow X$ of $\varphi$. 

 Define $\Psi:X\times \left[0,1\right]\rightarrow X$ by
\[
\Psi\left(x,t\right):=p\left(g_{1}^{-1}\Phi\left(g_{1}x,t\right),\ldots,g_{n}^{-1}\Phi\left(g_{n}x,t\right)\right).
\]
By Lemma~\ref{l:n mean +homotopy=G-homotopy}, 
$\Psi$ is an $H$-homotopy. Furthermore, since
 $\Phi(x,0)=\varphi(x,0)=x$, it follows from Lemma~\ref{l:n mean +homotopy=G-homotopy}-(1) that $\Psi_0$ is the identity map on $X$. 
 
 Now, if $x\in X^{H}$, then $g_{j}x\in X^{H}$ for all $g_{j}\in H$, so 
\[
\begin{split}
\Psi\left(x,t\right) & =p\left(g_{1}^{-1}\Phi\left(g_{1}x,t\right),\ldots,g_{n}^{-1}\Phi\left(g_{n}x,t\right)\right) \\
& =p\left(g_{1}^{-1}g_{1}x,\ldots,g_{n}^{-1}g_{n}x\right) \\
& =p\left(x,\ldots,x\right) \\
& =x.
\end{split}
\]

Then $\Psi\left(x,t\right)=x$ for all $(x,t)\in X^{H}\times[0,1]$.

Finally, let us prove that $\Psi_1$ is a retraction of $X$ onto $X^{H}$. Indeed, since $\Psi_1(x)=x$  for all $x\in X^{H}$, we only need to prove that $\Psi_{1}\left(X\right)\subseteq X^{H}$. Take $x\in X$ and notice that $hr\left(x\right)=r\left(x\right)$ for all $h\in H$, since $r\left(x\right)\in X^{H}$. This implies that 
\[
\begin{split}
h\Psi_{1}\left(x\right) & =hp\left(g_{1}^{-1}r\left(g_{1}x\right),\ldots,g_{n}^{-1}r\left(g_{n}x\right)\right) \\
& =p\left(hg_{1}^{-1}r\left(g_{1}x\right),\ldots,hg_{n}^{-1}r\left(g_{n}x\right)\right) \\
& =p\left(r\left(g_{1}x\right),\ldots,r\left(g_{n}x\right)\right) \\
& =p\left(g_{1}^{-1}r\left(g_{1}x\right),\ldots,g_{n}^{-1}r\left(g_{n}x\right)\right) \\
& =\Psi_{1}\left(x\right)
\end{split}
\]
for all $h\in H$.   Thus, $\Psi_{1}\left(x\right)\in X^{H}$, and therefore $\Psi_{1}$ is a retraction of $X$ onto $X^{H}$, as needed. We can now conclude that $X^{H}$ is an $H$-strong deformation retract of $X$.
\end{proof}

Regarding Lemma~\ref{l:fixed points are strong deformation retracts}-(1), it is interesting to point out that in general the set of fixed points may be empty, even if the group $G$ is finite and the space $X$ is an AR (and, therefore, contractible). 
\begin{example}
    Let $X:=\mathbb S_{\ell_2}$ be the unitary sphere of the Hilbert space $\ell_2$. By \cite[Chapter 6, Theorem 5.1]{bessaga}, $X$ is homeomorphic to $\ell_2$, and therefore it is an AR by Dugundji's Theorem (see, e.g., 
\cite[Theorem 1.5.1]{vanmill}). However, if we consider the continuous action of the group $\mathbb Z_2:=\{-1,1\}$ given by 
$$(t,x)\to tx,$$
then the set of $\mathbb Z_2$-fixed points is empty. In particular, $\mathbb S_{\ell_2}$ is not a $\mathbb Z_2$-AR.
\end{example}

\begin{theorem}\label{t:answer to Jaworowski}
Let $X$ be a metrizable $G$-space where $G$ is finite. Assume that for each subgroup $H\leq G$, the set $X^{H}$ is an AR. If there exists an equivariant $m$-mean for some multiple $m$ of $|G|$, then $X$ is a $G$-AR. 
\end{theorem}

\begin{proof}
Since $X=X^{\{e\}}$, we have that $X$ is an AR. Furthermore, by Lagrange's theorem, if $H\leq G$ is an arbitrary subgroup, $|H|$ divides $|G|$, and therefore $|H|$ divides $m$. Hence, according to Lemma~\ref{l: |G| divides n}, there exists an equivariant $|H|$-mean for every subgroup $H\leq G$. 

On the other hand,  since $X^{H}$ is closed in $X$ and $X^{H}$ is an $AR$,  there exists a retraction $r:X\to X^{H}$ for every subgroup $H\leq G$. Thus,  Lemma~\ref{l:fixed points are strong deformation retracts} implies that $X^{H}$ is an $H$-strong deformation retract of $X$ for every subgroup $H\leq G$. This, in combination with Theorem~\ref{t:Sergey strong deformation retract}, yields that $X$ is a $G$-AR.
\end{proof}

\begin{corollary}
Let $X$ be a metrizable $G$-space where $G$ is finite. Assume that for each subgroup $H$ of $G$, the set $X^{H}$ is an AR. If there exists an equivariant $|G|$-mean, then $X$ is a $G$-AR. 
\end{corollary}

Regarding Question~\ref{preguntahugo} (\cite[Question 1]{hugo}), in the following theorem and corollary, we show that we can obtain a positive answer if the group is finite and $n$ is a multiple of $|G|$.

\begin{theorem} \label{t:n mean implica G-ANR} Let $G$ be a finite group and let $X$ be 
 a connected $G$-ANR with finitely generated homology groups such that almost all vanish. Assume that there exists an equivariant $n$-mean $p:X^{n}\rightarrow X$ for some multiple $n$ of $|G|$, then $X$  is a $G$-AR.
\end{theorem}

\begin{proof}
By Lemma~\ref{l: |G| divides n}, we can assume without loss of generality, that $n=|G|$. 

Since $X$ is a $G$-ANR, it is an ANR (\cite[Theorem 3.7]{s3}). Then, it follows from Theorem~\ref{t:caracterizacion n mean original}\big((2)$\Rightarrow$(3)\big) that $X$ is an AR and therefore is contractible. 

Let $\Phi:X\times \left[0,1\right]\rightarrow X$ be a contraction of $X$, and let $g_{1},\ldots,g_{n}$ be the elements of $G$. 
Consider the function $\Psi:X\times \left[0,1\right]\rightarrow X$ given by
\[
\Psi\left(x,t\right)=p\left(g_{1}^{-1}\Phi\left(g_{1}x,t\right),\ldots,g_{n}^{-1}\Phi\left(g_{n}x,t\right)\right).
\]
By Lemma~\ref{l:n mean +homotopy=G-homotopy}, $\Psi$ is an equivariant contraction and therefore $X$ is $G$-contractible. 
Then,  by Theorem \ref{gc}, $X$ is a $G$-AR. 
\end{proof}

If $X$ is a compact connected $G$-ANR,  it has finitely generated homology groups and almost all vanish (\cite[Chapter 4, Corollary 7.2]{Hu}). Then, as a corollary of Theorem~\ref{t:n mean implica G-ANR}, we obtain the following. 

\begin{corollary}\label{c:problema hugo}
 Let $G$ be a finite group, and let $X$ be 
 a compact connected $G$-ANR.  If there exists an equivariant $n$-mean $p:X^{n}\rightarrow X$ for some multiple $n$ of $|G|$ (in particular, if $n=|G|$), then $X$ is a $G$-AR. 
\end{corollary}

\section{Quasi-means}

We start this section by noticing that in Theorem~\ref{t:caracterizacion n mean original}, none of the conditions (M1) or (M2) can be dropped. 

\begin{example}\label{e:circulo}
Consider the circle $\mathbb{S}^1:=\{(x,y)\in\mathbb R^{2}: x^2+y^2=1\}$, which is a compact connected ANR that is not contractible. Let $z_0\in\mathbb S^{1}$ be fixed. Then, for every $n\in\mathbb N$, the constant map $c:\mathbb S^{1}\times\mathbb S^{1}\to\mathbb S^{1}$ given by $c(x,y)=z_0$ satisfies condition (M2) (but not (M1)).

On the other hand, the dictatorial map $p:\mathbb{S}^1\times \mathbb{S}^1\rightarrow \mathbb{S}^1$ defined by $p\left(x,y\right)=x$ is not a mean (because it does not satisfy condition (M2)), but it satisfies condition (M1) and therefore is a quasi-mean.
\end{example}

After this example, it is natural to ask what modifications can we  make to the definition of an $n$-mean, so that the existence of such a map guarantees the contractibility of the space. 

Since we want to avoid solomonic maps, we will be interested in  continuous maps $p:X^{n}\rightarrow X$ satisfying the following inequality for a fixed  $\lambda\in(0,1)$:

\begin{equation}\label{eq: contractive quasi mean}
\max_{i=1,\ldots,n}d\left(x_{i},p\left(x_{1},\ldots,x_{n}\right)\right)\leq \lambda \max_{j,k=1,\ldots,n}d\left(x_{j},x_{k}\right) 
\end{equation}
where $\left(x_{1},\ldots,x_{n}\right)\in X^{n}$.
%then Question \ref{preguntahugo} has a positive answer. In Theorem \ref{cont} we prove the case when $n=2$, and in Corollary \ref{contc} we see the general case.

%\begin{example}
%Let $G$ be a compact topological group acting on the circle $\mathbb{S}^1:=\{(x,y)\in\mathbb R^{2}: x^2+y^2=1\}$. Consider the dictatorial map $p:\mathbb{S}^1\times \mathbb{S}^1\rightarrow \mathbb{S}^1$ defined by $p\left(x,y\right)=x$. This is an equivariant quasi-mean that satisfies that 
% \[
%\max\left\{d\left(x,p\left(x,y\right)\right),d\left(y,p\left(x,y\right)\right)\right\}=d\left(x,y\right)
%\]
%for every $x,y\in \mathbb{S}^1$. However, $\mathbb{S}^1$ is not $G$-contractible because it is not even contractible. 
%\end{example}

Notice that if a continuous function $p:X^{n}\rightarrow X$ satisfies inequality~\ref{eq: contractive quasi mean}, 
%\[
%\max_{i=1,\ldots,n}d\left(x_{i},p\left(x_{1},\ldots,x_{n}\right)\right)\leq \lambda \max_{j,k=1,\ldots,n}d\left(x_{j},x_{k}\right) 
%\]
for some positive $\lambda$, then 
$$0\leq d(x,p(x,\dots,x))\leq \lambda d(x,x)=0.$$
This implies that $p$ satisfies condition (M1) and therefore 
it is an $n$-quasi-mean. 
This situation motivates the following definition.

\begin{definition}\label{d:quasimean}
    Let $(X,d)$ be a metric space. A continuous function $p:X^{n}\rightarrow X$ is a \textit{contractive $n$-quasi-mean} if there exists $\lambda\in \left(0,1\right)$ such that 
    \[
\max_{i=1,\ldots,n}d\left(x_{i},p\left(x_{1},\ldots,x_{n}\right)\right)\leq \lambda \max_{j,k=1,\ldots,n}d\left(x_{j},x_{k}\right) 
\]
for every $\left(x_{1},\ldots,x_{n}\right)\in X^{n}$. If $n=2$, we simply say that $p$ is a \textit{contractive quasi-mean}.
\end{definition}

At first glance, inequality~(\ref{eq: contractive quasi mean}) may seem a little artificial. However, it captures a desirable property: that the output of the values $x_1,\dots ,x_n$ lies \textit{in between} those values. Furthermore, this inequality holds in some classic situations. For instance, if $X$ is a convex subset of a normed space, then it is not difficult to prove that the arithmetic $n$-mean
$A:X^n\to X$ given by 
$$A(x_1,\dots , x_n):=\frac{1}{n}\sum_{i=1}^nx_i$$
satisfies inequality~(\ref{eq: contractive quasi mean}) for $\lambda=\frac{n-1}{n}$.

Also, if $I=\left[a,b\right]$ is an interval contained in the set of positive real numbers $\mathbb{R}^{+}$, then the geometric mean $p:I\times I\rightarrow I$, given by $p(x,y)=\sqrt{xy}$, satisfies inequality~(\ref{eq: contractive quasi mean}) for $\lambda=\frac{b-\sqrt{ab}}{b-a}$.

Throughout the rest of the section, we will use the following notation.  For every $n\in \mathbb{N}$, we will denote by $D_{n}$ the set
$$D_n:=\left\{\frac{j}{2^{n}} \mid j\in \left\{0,1,\ldots,2^{n}\right\}\right\},$$
and let $D:=\bigcup_{n=0}^{\infty}D_{n}$ be the set of all dyadic rationals in $\left[0,1\right]$. Furthermore,
for every $x\in D$ define 
\begin{equation}\label{eq: h diadicos}
    h\left(x\right):=\min\left\{n\in \mathbb{N} \mid x\in D_{n}\right\}.
\end{equation}

\begin{lemma} \label{st}
Suppose that $s,t\in D$ and $s<t$. If $h$ is the map defined in equation~(\ref{eq: h diadicos}), then there exist $k,l\in \mathbb{N}$ and $s_{0},s_{1},\ldots,s_{k},t_{0},t_{1},\ldots,t_{l}\in D$ satisfying the following conditions:
\begin{enumerate}[\rm(1)]
\item $s=s_{0}\leq\ldots\leq s_{k}=t_{l}\leq \ldots\leq t_{0}=t$.
\item $h\left(s_{0}\right)>\ldots>h\left(s_{k}\right)$ and $h\left(t_{0}\right)>\ldots>h\left(t_{l}\right)$.
\item For every $i\in \left\{0,\ldots,k-1\right\}$, $s_{i},s_{i+1}\in D_{h\left(s_{i}\right)}$ and $\left|s_{i}-s_{i+1}\right|=\frac{1}{2^{h\left(s_{i}\right)}}$.
\item For every $i\in \left\{0,\ldots,l-1\right\}$, $t_{i},t_{i+1}\in D_{h\left(t_{i}\right)}$ and $\left|t_{i}-t_{i+1}\right|=\frac{1}{2^{h\left(t_{i}\right)}}$.
\end{enumerate}
\end{lemma}

\begin{proof}
For every $s,t\in D$ we have that $s,t\in D_{\max\left\{h\left(s\right),h\left(t\right)\right\}}$, so there is a unique $i\in \mathbb{N}$ such that $\left|s-t\right|=\frac{i}{2^{\max\left\{h\left(s\right),h\left(t\right)\right\}}}$. Let us denote $\rho\left(s,t\right):=i$. We will prove the result by induction on $\rho\left(s,t\right)$.

Let $s,t\in D$ with $s<t$ and suppose that $\rho\left(s,t\right)=1$. Then there exists $j\in \mathbb{N}$ such that $s=\frac{j}{2^{\max\left\{h\left(s\right),h\left(t\right)\right\}}}$ and $t=\frac{j+1}{2^{\max\left\{h\left(s\right),h\left(t\right)\right\}}}$. Notice that, since either $j$ or $j+1$ is even, then $h\left(s\right)\neq h\left(t\right)$. If $h\left(s\right)>h\left(t\right)$, then $s_{0}:=s$ and $s_{1}:=t=:t_{0}$ satisfy the properties (1)-(4). On the other hand,  if $h\left(t\right)>h\left(s\right)$, then $s_{0}:=s=:t_{1}$ and $t_{0}:=t$ are the required elements of $D$.

Let $n\geq 1$ and suppose that the result holds whenever $\rho\left(s,t\right)\leq n$. Now, take $s,t\in D$ such that $s<t$ and $\rho\left(s,t\right)=n+1$. Define $s_{0}:=s$ and $t_{0}:=t$. 

If $h\left(s_{0}\right)\geq h\left(t_{0}\right)$, we  define $s_{1}:=s_{0}+\frac{1}{2^{h\left(s_{0}\right)}}$. Since $s_{0},t_{0}\in D_{h\left(s_{0}\right)}$, then $\left|t_{0}-s_{0}\right|=t_0-s_0\geq \frac{1}{2^{h\left(s_{0}\right)}}$, which implies that $s_{1}\leq t_{0}$. Also, notice that $s_{0}=\frac{m}{2^{h\left(s_{0}\right)}}$ for some odd $m\in \mathbb{N}$, so $s_{1}=\frac{m+1}{2^{h\left(s_{0}\right)}}$ with $m+1$ even. This implies that $h\left(s_{0}\right)>h\left(s_{1}\right)$. So, if $s_{1}=t_{0}$, then $s_{0},s_{1},t_{0}$ satisfy properties (1)-(4). 
If $s_{1}<t_{0}$,  since $\rho\left(s_{1},t_{0}\right)< \rho\left(s_{0},t_{0}\right)$, we can use our induction hypothesis and find $s_{1}=s_{0}^{\prime},\ldots,s_{k}^{\prime},t_0=t_{0}^{\prime},\ldots,t_{l}^{\prime}\in D$ satisfying properties (1)-(4) for   $s_{1}$ and $t_{0}$. Therefore, $s_{0},s_{0}^{\prime},\ldots,s_{k}^{\prime},t_{0}^{\prime},\ldots,t_{l}^{\prime}$ are the required elements of $D$.

If $h\left(s_{0}\right)<h\left(t_{0}\right)$, we define $t_{1}:=t_{0}-\frac{1}{2^{h\left(t_{0}\right)}}$ and then use an argument analogous to the one given in the previous paragraph.
\end{proof}

\begin{theorem} \label{cont}
Let $\left(X,d\right)$ be a complete metric space. 
\begin{enumerate}
    \item If there exists a contractive quasi-mean $p:X\times X\rightarrow X$, then $X$ is contractible.
    \item If, in addition, $X$ is a $G$-space such that $X^{G}$ is non-empty and $p$ is equivariant, then $X$ is $G$-contractible.
\end{enumerate}
\end{theorem}

\begin{proof}
(1) Let $\lambda\in \left(0,1\right)$ be such that
\[
 \max\left\{d\big(x,p(x,y)\big),d\big(y,p(x,y)\big)\right\}\leq \lambda d\left(x,y\right)
    \]
for every $x,y\in X$.

Fix $\theta\in X$ and consider $\phi_{0}:X\times D_{0}\rightarrow X$ given by $\phi_{0}\left(x,0\right)=x$ and $\phi_{0}\left(x,1\right)=\theta$. For every $n\geq 1$, define 
recursively $\phi_{n}:X\times D_{n}\rightarrow X$ by
\[
\phi_{n}\left(x,\frac{j}{2^{n}}\right)=
\left\{\begin{array}{ll} 
\phi_{n-1}\left(x,\frac{j}{2^{n}}\right) & \mbox{ if $j$ is even},\\
p\left(\phi_{n-1}\left(x,\frac{j-1}{2^{n}}\right),\phi_{n-1}\left(x,\frac{j+1}{2^{n}}\right)\right) & \mbox{ if $j$ is odd}.\end{array}
\right. 
\]
It can be checked by induction that $\phi_{n}$ is continuous for every $n\in \mathbb{N}$. Now we can define $\Phi:X\times D\rightarrow X$ by $\Phi\left(x,\frac{j}{2^{n}}\right)=\phi_{n}\left(x,\frac{j}{2^{n}}\right)$. Then $\Phi\restriction_{X\times\left\{t\right\}}$ is continuous for each $t\in D$.

For every $x\in X$, let $\Phi_{x}:D\rightarrow X$ be given by $\Phi_{x}\left(t\right):=\Phi\left(x,t\right)$. We will prove that each $\Phi_{x}$ is uniformly continuous. Take $x\in X$. 

\textbf{Claim 1.} For every $n\in \mathbb{N}$ and $j\in \left\{0,1,\ldots,2^{n}-1\right\}$,

\[
d\left(\Phi_{x}\Bigl(\frac{j}{2^{n}}\Bigr),\Phi_{x}\Bigl(\frac{j+1}{2^{n}}\Bigr)\right)\leq \lambda^{n}d\left(x,\theta\right).
\]

\textit{Proof of Claim 1.}
We  proceed by induction on $n$. Since $$d\big(\Phi_{x}(0),\Phi_{x}(1)\big)=d\left(x,\theta\right),$$ the claim is true for $n=0$.

Let $m\geq 1$. Assume the claim is true for $m-1$ and let $j\in \left\{0,1,\ldots,2^{m}-1\right\}$. We will suppose that $j$ is odd, since the case where $j$ is even is analogous. Using the fact that $p$ is a contractive quasi-mean and our induction hypothesis, we obtain that 
\[
\begin{split}
d\left(\Phi_{x}\Bigl(\frac{j}{2^{m}}\Bigr),\Phi_{x}\Bigl(\frac{j+1}{2^{m}}\Bigr)\right) & =d\left(p\left(\Phi_{x}\Bigl(\frac{j-1}{2^{m}}\Bigr),\Phi_{x}\Bigl(\frac{j+1}{2^{m}}\Bigr)\right),\Phi_{x}\Bigl(\frac{j+1}{2^{m}}\Bigr)\right) \\
& \leq \lambda d\left(\Phi_{x}\Bigl(\frac{j-1}{2^{m}}\Bigr),\Phi_{x}\Bigl(\frac{j+1}{2^{m}}\Bigr)\right) \\
&=\lambda d\left(\Phi_{x}\biggl(\frac{\frac{j-1}{2}}{2^{m-1}}\biggr),\Phi_{x}\biggl(\frac{\frac{j+1}{2}}{2^{m-1}}\biggr)\right)\\
&\leq \lambda\lambda^{m-1}d\left(x,\theta\right)= \lambda^{m}d\left(x,\theta\right).
\end{split}
\]
Therefore, the claim is true for every $n\in \mathbb{N}$.

\QEDA

\textbf{Claim 2.} For every $s,t\in D$,
\begin{equation} \label{claim2}
    d\big(\Phi_{x}(s),\Phi_{x}(t)\big)\leq \frac{2 d\left(x,\theta\right)}{1-\lambda}\left|s-t\right|^{\alpha},
\end{equation}
where $\alpha:=-\frac{\ln\left(\lambda\right)}{\ln\left(2\right)}>0$.

\textit{Proof of Claim 2.}
Take $s,t\in D$ with $s<t$, and $s=s_{0},\ldots,s_{k},t=t_{0},\ldots,t_{l}\in D$ satisfying the properties (1)-(4) of Lemma \ref{st}. Notice that, by Claim 1 and properties (3) and (4) of Lemma \ref{st}, 
\[
d\big(\Phi_{x}(s_{i}),\Phi_{x}(s_{i+1})\big)\leq \lambda^{h\left(s_{i}\right)}d\left(x,\theta\right)
\]
for every $i\in \left\{0,\ldots,k-1\right\}$ and
\[
d\big(\Phi_{x}(t_{i}),\Phi_{x}(t_{i+1})\big)\leq \lambda^{h\left(t_{i}\right)}d\left(x,\theta\right)
\]
for every $i\in \left\{0,\ldots,l-1\right\}$. Now, by Lemma \ref{st}, $h(s_{0})>\ldots>h(s_{k-1})$ and $h(t_{0})>\ldots>h(t_{l-1})$, which implies that 
\begin{align}
d\big(\Phi_{x}(s),\Phi_{x}(t)\big) & \leq \sum_{i=0}^{k-1}d\big(\Phi_{x}(s_{i}),\Phi_{x}(s_{i+1})\big)+\sum_{i=0}^{l-1}d\big(\Phi_{x}(t_{i}),\Phi_{x}(t_{i+1})\big) \nonumber\\
& \leq \sum_{i=0}^{k-1}\lambda^{h\left(s_{i}\right)}d\left(x,\theta\right)+\sum_{i=0}^{l-1}\lambda^{h\left(t_{i}\right)}d\left(x,\theta\right) \nonumber\\
& \leq d\left(x,\theta\right) \left(\sum_{n=h\left(s_{k-1}\right)}^{\infty}\lambda^{n}+\sum_{n=h\left(t_{l-1}\right)}^{\infty}\lambda^{n}\right) \nonumber\\
& =d\left(x,\theta\right)\frac{1}{1-\lambda}\left(\lambda^{h\left(s_{k-1}\right)}+\lambda^{h\left(t_{l-1}\right)}\right). \label{des: distancia phix-1}
\end{align}

 Also notice that $\frac{1}{2^{h\left(s_{k-1}\right)}}=\left|s_{k-1}-s_{k}\right|\leq \left|s-t\right|$, therefore
\[
\lambda^{h\left(s_{k-1}\right)}=e^{\ln\left(\lambda\right)h\left(s_{k-1}\right)}=e^{-\alpha\ln\left(2\right)h\left(s_{k-1}\right)}=\left(\frac{1}{2^{h\left(s_{k-1}\right)}}\right)^{\alpha}\leq \left|s-t\right|^{\alpha},
\]
where $\alpha=-\frac{\ln\left(\lambda\right)}{\ln\left(2\right)}$.
Using the same argument, we can prove that $\lambda^{h\left(t_{l-1}\right)}\leq \left|s-t\right|^{\alpha}$. Thus, we can conclude that
\[
d\big(\Phi_{x}(s),\Phi_{x}(t)\big)\leq \frac{2 d\left(x,\theta\right)}{1-\lambda}\left|s-t\right|^{\alpha},
\]
and this finishes the proof of the claim.
\QEDA

Claim 2 implies that $\Phi_{x}:D\rightarrow X$ is Hölder continuous and therefore uniformly continuous.
%because, since $r\rightarrow \frac{2 d\left(x,\theta\right)}{1-\lambda}r^{\alpha}$ is a continuous function, particularly at $0$, for every $\varepsilon>0$ there exists $\delta>0$ such that if $\left|s-t\right|<\delta$ then $\frac{2 d\left(x,\theta\right)}{1-\lambda}\left|s-t\right|^{\alpha}<\varepsilon$.
Since $D$ is dense in $\left[0,1\right]$ and $(X,d)$ is complete, there exists  a continuous extension $\widetilde{\Phi}_{x}:\left[0,1\right]\rightarrow X$ of $\Phi_{x}$ (\cite[Theorem 4.3.17]{engelking}). Notice that, by taking limits on both sides of inequality (\ref{claim2}), we can  prove that $\widetilde{\Phi}_{x}$ satisfies the same inequality. That is, for every $s,t\in \left[0,1\right]$,
\begin{equation} \label{3}
d\big(\widetilde{\Phi}_{x}(s),\widetilde{\Phi}_{x}(t)\big)\leq \frac{2 d\left(x,\theta\right)}{1-\lambda}\left|s-t\right|^{\alpha}.
\end{equation}

Let us define $\widetilde{\Phi}:X\times \left[0,1\right]\rightarrow X$ by $\widetilde{\Phi}\left(x,t\right)=\widetilde{\Phi}_{x}\left(t\right)$. Notice that, for every $t\in D$, $\widetilde{\Phi}\left(x,t\right)=\widetilde{\Phi}_{x}\left(t\right)=\Phi_{x}\left(t\right)=\Phi\left(x,t\right)$, so $\widetilde{\Phi}$ is an extension of $\Phi$. This means that $\widetilde{\Phi}\left(x,0\right)=x$ and $\widetilde{\Phi}\left(x,1\right)=\theta$ for every $x\in X$. Hence, in order to prove that $\widetilde{\Phi}$ is a homotopy, and  therefore that $X$ is contractible, it only remains to show that $\widetilde{\Phi}$ is continuous. To prove this, take a point $(x,\eta)\in X\times \left[0,1\right]$ and suppose that $\left((x_{n},t_{n})\right)_{n\in\mathbb N}\subset X\times \left[0,1\right] $ is a sequence that converges to $\left(x,\eta\right)$ in $X\times \left[0,1\right]$. We shall prove that $\big(\widetilde{\Phi}(x_n, t_n)\big)_{n\in\mathbb N}$ converges to $\widetilde{\Phi}(x,\eta)$.

\textbf{Claim 3.} The family $\mathcal{H}:=\left\{\widetilde{\Phi}_{x_{n}} \mid n\in \mathbb{N}\right\}\cup \left\{\widetilde{\Phi}_{x}\right\}$ is uniformly equicontinuous.

\textit{Proof of Claim 3.}
Since $\left(x_{n}\right)_{n\in\mathbb N}$ converges to $x$, the set 
\begin{equation}\label{eq:set A}
    A:=\left\{x_{n} \mid n\in \mathbb{N}\right\}\cup\{x\}
\end{equation} is compact and therefore we can find $M>0$ such that $A\subseteq B\left(\theta,M\right)$. Then, by inequality (\ref{3}), for every $z\in A$ and every pair $s,t\in [0,1]$, we have that 
\begin{equation}\label{eq: familia uniformemente acotada}
d\big(\widetilde{\Phi}_{z}(s),\widetilde{\Phi}_{z}(t)\big)\leq \frac{2 M}{1-\lambda}\left|s-t\right|^{\alpha}.
\end{equation}

Let $\varepsilon>0$ and pick $\delta>0$ such that,   if $s,t\in[0,1]$ and $\left|s-t\right|<\delta$, then $\frac{2M}{1-\lambda}\left|s-t\right|^{\alpha}<\varepsilon$. This, in combination with inequality~(\ref{eq: familia uniformemente acotada}), guarantees that  $d\big(\widetilde{\Phi}_{z}(s),\widetilde{\Phi}_{z}(t)\big)<\varepsilon$ for every $z\in A$ and every $s,t\in [0,1]$ with $|s-t|<\delta$. This means that $\mathcal{H}$ is a uniformly equicontinuous family.

\QEDA

\textbf{Claim 4.} The sequence $\big(\widetilde{\Phi}_{x_{n}}(t)\big)$ converges to $\widetilde{\Phi}_{x}\left(t\right)$ for every $t\in \left[0,1\right]$.

\textit{Proof of Claim 4.}
If $t\in D$, we know that $\Phi\restriction_{X\times\left\{t\right\}}$ is continuous. Thus,  the sequence $\big(\Phi(x_{n},t)\big)_{n\in\mathbb N}=\big(\widetilde{\Phi}_{x_{n}}(t)\big)_{n\in\mathbb N}$ converges to $\Phi(x,t)=\widetilde{\Phi}_{x}(t)$.

Now, take $t\in \left[0,1\right]\backslash D$ and suppose that $\big(\widetilde{\Phi}_{x_{n}}(t)\big)_{n\in\mathbb N}$ does not converge to $\widetilde{\Phi}_{x}(t)$. Then there exist $\varepsilon_0>0$ and a subsequence $\left(x_{n_{l}}\right)_{l\in\mathbb N}$ of $\left(x_{n}\right)_{n\in\mathbb N}$ such that $d\big(\widetilde{\Phi}_{x_{n_{l}}}(t),\widetilde{\Phi}_{x}(t)\big)\geq \varepsilon_0$ for every $l\in \mathbb{N}$. 
By Claim 3, we can find $\delta>0$ such that, if $s\in \left[0,1\right]$ and $\left|t-s\right|<\delta$, then $d\big(\widetilde{\Phi}_{z}(t),\widetilde{\Phi}_{z}(s)\big)< \frac{\varepsilon_0}{3}$ for every $z\in A$ (where $A$ is the set defined in equation~(\ref{eq:set A})). Take $r\in D\cap (t-\delta,t+\delta)$ and observe that 
\[
\begin{split}
\varepsilon_0 & \leq d\big(\widetilde{\Phi}_{x_{n_{l}}}(t),\widetilde{\Phi}_{x}(t)\big) \\ 
& \leq d\big(\widetilde{\Phi}_{x_{n_{l}}}(t),\widetilde{\Phi}_{x_{n_{l}}}(r)\big)+d\big(\widetilde{\Phi}_{x_{n_{l}}}(r),\widetilde{\Phi}_{x}(r)\big)+d\big(\widetilde{\Phi}_{x}(r),\widetilde{\Phi}_{x}(t)\big) \\ 
& < d\big(\widetilde{\Phi}_{x_{n_{l}}}(r),\widetilde{\Phi}_{x}(r)\big)+\frac{2\varepsilon_0}{3},
\end{split}
\]
for every $l\in \mathbb N$.
Therefore, $d\big(\widetilde{\Phi}_{x_{n_{l}}}(r),\widetilde{\Phi}_{x}(r)\big)\geq \frac{\varepsilon_0}{3}$ for every $l\in \mathbb{N}$, but this contradicts the fact that $\big(\widetilde{\Phi}_{x_{n}}(r)\big)$ converges to $\widetilde{\Phi}_{x}\left(r\right)$. Then we can conclude that $\big(\widetilde{\Phi}_{x_{n}}(t)\big)$ converges to $\widetilde{\Phi}_{x}\left(t\right)$, as desired. 
\QEDA

After Claim 4, we conclude that $\big(\widetilde{\Phi}_{x_n}\big)_{n\in\mathbb N}$ converges pointwise to $\widetilde{\Phi}_x$. However, since $\mathcal{H}$ is an equicontinuous family and $[0,1]$ is compact, this convergence is also uniform (see, e.g., \cite[Chapter XII, Problem Section 7-1]{dugundji} or \cite[Theorem 43.14]{willard}).
Hence, by \cite[Theorem 7.5, Chapter XII]{dugundji}, the sequence $\big(\widetilde{\Phi}_{x_n}(t_n)\big)_{n\in\mathbb N}=\big(\widetilde{\Phi}(x_n, t_n)\big)_{n\in\mathbb N}$ converges to $\widetilde{\Phi}(x,\eta)$, as needed. 
%By Claim 3, $\mathcal{H}=\left\{\widetilde{\Phi}_{x_{n}} \mid n\in \mathbb{N}\right\}\cup \left\{\widetilde{\Phi}_{x}\right\}$ is an equicontinuous family, and by Claim 4, $\left\{\widetilde{\Phi}_{x_{n}}\left(t\right) \mid n\in\mathbb{N}\right\}\cup \left\{\widetilde{\Phi}_{x}\left(t\right)\right\}$ is compact for every $t\in \left[0,1\right]$. Therefore, by the Arzelà-Ascoli theorem (\cite[Chapter 12, Theorem 6.4]{dugundji}), $\mathcal{H}$ is precompact in C$\left(\left[0,1\right],X\right)$, the space of continuous functions from $\left[0,1\right]$ to $X$ equipped with the compact-open topology. Then, there exists a subsequence $\left(\widetilde{\Phi}_{x_{n_{k}}}\right)$ of $\left(\widetilde{\Phi}_{x_{n}}\right)$ that converges uniformly and, by Claim 4, it must converge to $\widetilde{\Phi}_{x}$. This implies that $\left(\widetilde{\Phi}_{x_{n_{k}}}\left(t_{n_{k}}\right)\right)$ converges to $\widetilde{\Phi}_{x}\left(t\right)$. That is, $\left(\widetilde{\Phi}\left(x_{n_{k}},t_{n_{k}}\right)\right)$ converges to $\widetilde{\Phi}\left(x,t\right)$.

%We have proved that, for any sequence $\left(\left(x_{n},t_{n}\right)\right)$ that converges to some $\left(x,t\right)$ in $X\times \left[0,1\right]$, there exists a subsequence $\left(\left(x_{n_{k}},t_{n_{k}}\right)\right)$ of $\left(\left(x_{n},t_{n}\right)\right)$ such that $\left(\widetilde{\Phi}\left(x_{n_{k}},t_{n_{k}}\right)\right)$ converges to $\widetilde{\Phi}\left(x,t\right)$. We can conclude that $\widetilde{\Phi}$ is continuous. Therefore, $X$ is contractible.

(2) Now, suppose  that $X$ is a $G$-space such that $X^{G}$ is non-empty and that $p$ is equivariant. Since $X^{G}$ is non-empty, we may assume that $\theta$, the point we used to define the homotopy $\widetilde{\Phi}$, is a $G$-fixed point. This implies that  $\phi_{0}\restriction_{X\times\left\{t\right\}}$ is equivariant for $t=0$ and $t=1$, and it can be checked by induction that $\phi_{n}\restriction_{X\times\left\{t\right\}}$ is equivariant for each $t\in D_{n}$. Therefore, $\Phi\restriction_{X\times\left\{t\right\}}$ is equivariant.

To prove that $X$ is $G$-contractible, it remains to show that $\widetilde{\Phi}\restriction_{X\times\left\{t\right\}}$ is equivariant for every $t\in \left[0,1\right]$. Take $t\in \left[0,1\right]$, $x\in X$, and $g\in G$. Let $\left(r_{n}\right)$ be a sequence in $D$ that converges to $t$. Since $\Phi\restriction_{X\times\left\{r_{n}\right\}}$ is equivariant for each $n\in \mathbb{N}$, then
\begin{align*}
\widetilde{\Phi}\left(gx,t\right)&=\widetilde{\Phi}\left(gx,\lim_{n\to\infty}r_{n}\right)=\lim_{n\to\infty}\widetilde{\Phi}\left(gx,r_{n}\right)=\lim_{n\to\infty}\Phi\left(gx,r_{n}\right)\\
&=\lim_{n\to\infty}g\Phi\left(x,r_{n}\right)=g\lim_{n\to\infty}\widetilde{\Phi}\left(x,r_{n}\right)=g\widetilde{\Phi}\left(x,t\right).
\end{align*}
Therefore, $\widetilde{\Phi}\restriction_{X\times\left\{t\right\}}$ is equivariant.
\end{proof}

\begin{corollary} \label{contc}
Let $\left(X,d\right)$ be a complete metric space. 
\begin{enumerate}
    \item If for some $n\in \mathbb{N}$ there exists a contractive $n$-quasi-mean  $p:X^{n}\rightarrow X$, then $X$ is contractible.
    \item If, in addition, $X$ is a $G$-space such that $X^{G}$ is non-empty and $p$ is equivariant, then $X$ is $G$-contractible. Moreover, if $X$ is a $G$-ANR, then it is a $G$-AR.
\end{enumerate}
\end{corollary}

\begin{proof}
Define $p^{\prime}:X\times X\rightarrow X$ by $p^{\prime}\left(x,y\right)=p\left(x,y,\ldots,y\right)$. Then, $p^{\prime}$ is a contractive quasi-mean. Therefore, by Theorem \ref{cont}, $X$ is contractible. 

 If, in addition, $X$ is a $G$-space such that $X^{G}$ is non-empty and $p$ is equivariant, then $p^{\prime}$ is also equivariant and, by Theorem \ref{cont}, $X$ is $G$-contractible.  Hence, if $X$ is also a $G$-ANR, Theorem~\ref{gc} implies that $X$ is a $G$-AR, as desired. 
\end{proof}

We finish this work by noticing that the inequality~\ref{eq: contractive quasi mean} cannot be weakened to
\[
\max_{i=1,\ldots,n}d\left(x_{i},p\left(x_{1},\ldots,x_{n}\right)\right)\leq\max_{j,k=1,\ldots,n}d\left(x_{j},x_{k}\right) 
\]
for every $\left(x_{1},\ldots,x_{n}\right)\in X^{n}$. Indeed, the dictatorial map $p:\mathbb{S}^1\times \mathbb{S}^1\rightarrow \mathbb{S}^1$ defined in  Example~\ref{e:circulo}, satisfies 
 \[
\max\left\{d\left(w,p\left(w,v\right)\right),d\left(v,p\left(w,v\right)\right)\right\}=d\left(w,v\right),
\]
where $d$ stands for the Euclidean distance in $\mathbb{S}^1$. However, $\mathbb{S}^1$ is not contractible.

It is natural to ask whether $X$ being compact and $p:X^{n}\rightarrow X$ satisfying
 \begin{equation} \label{eq:nolambda}
\max_{i=1,\ldots,n}d\left(x_{i},p\left(x_{1},\ldots,x_{n}\right)\right)<\max_{j,k=1,\ldots,n}d\left(x_{j},x_{k}\right) 
 \end{equation}
for every $\left(x_{1},\ldots,x_{n}\right)\in X^{n}$ imply that $p$ is a contractive $n$-quasi-mean. For instance, the geometric mean $p:\mathbb{R}^{+}\times \mathbb{R}^{+}\rightarrow \mathbb{R}^{+}$, given by $p(x,y)=\sqrt{xy}$, satisfies (\ref{eq:nolambda}), but not (\ref{eq: contractive quasi mean}) for any $\lambda\in \left(0,1\right)$. However, as we have already noticed,  if we restrict $p$ to a compact interval contained in $\mathbb{R}^{+}$, then it satisfies (\ref{eq: contractive quasi mean}). But this is not true for every compact space, as the following example shows. 

\begin{example}
Consider $p:\left[0,1\right]\times \left[0,1\right]\rightarrow \left[0,1\right]$ given by 
\[
p\left(x,y\right)=\min\left\{x,y\right\}+\frac{|x-y|^{2}}{2}.
\]
\end{example}
Then $|x-p\left(x,y\right)|<|x-y|$ for every $x,y\in \left[0,1\right]$, thus $p$ satisfies (\ref{eq:nolambda}). However, notice that
\[
\frac{|\frac{1}{n}-p\left(\frac{1}{n},0\right)|}{|\frac{1}{n}-0|}=1-\frac{1}{2n}\to 1
\]
as $n\to \infty$, and therefore there is no $\lambda\in \left(0,1\right)$ such that $|x-p\left(x,y\right)|\leq\lambda|x-y|$ for every $x,y\in \left[0,1\right]$. Then, $p$ is not a contractive quasi-mean.


\begin{thebibliography}{23}
\bibitem{s5}
S. Antonyan. \textit{Retracts in categories of $G$-spaces}. Sov. J. Contemp. Math. Anal. 15 (1980), 365-378.
\bibitem{s1}
S. Antonyan. \textit{Equivariant embeddings into G-AR's}. Glas. Mat. 22 (1987), 503-533.
\bibitem{s4}
S. Antonyan. \textit{On based-free actions of compact Lie groups on the Hilbert cube}. Math. Notes 65, 2 (1999), 135-143. \url{https://doi.org/10.1007/BF02679809}
\bibitem{s3}
S. Antonyan. \textit{A characterization of equivariant absolute extensors and the equivariant Dugundji Theorem}. Houst. J. Math. 31 (2005), 451-462.
\bibitem {s2}
S. Antonyan. \textit{Some open problems in equivariant infinite-dimensional topology}. Topol. Appl. 311 (2022). \url{https://doi.org/10.1016/j.topol.2021.107966}

\bibitem{Arrow} K. J. Arrow, \textit{Social Choice and Individual Values}, Yale Univ. Press. New Haven,
Conn., 1951; 2nd ed., Wiley, New York, 1963.
\bibitem {bessaga}
C. Bessaga and A. Pelczynski. \textit{Selected Topics in Infinite-Dimensional Topology}. Polish Scientific Publishers, Warsaw, 1975.
\bibitem{Bredon} G. E. Bredon. \textit{Introduction to Compact Transformation Groups}. Academic Press,
New York, 1972.
\bibitem{ch}
G. Chichilnisky and G. Heal. \textit{Necessary and sufficient conditions for a resolution of the social choice paradox}. J. Econ. Theory 31 (1983), 68-87. \url{https://doi.org/10.1016/0022-0531(83)90021-2}
\bibitem{dugundji}
J. Dugundji. \textit{Topology}. Allyn and Bacon, Boston, 1966.
\bibitem{eck}
B. Eckmann. \textit{Räume mit Mittelbildungen}. Comment. Math. Helv. 28 (1954), 329-340.
\bibitem{eck2}
B. Eckmann, T. Ganea, P. J. Hilton. \textit{Generalized means}. Studies in Mathematical Analysis, Stanford University Press (1962), 82-92.
\bibitem{eck3} B. Eckmann. \textit{Social choice and topology. A case of pure and applied mathematics}. Expo. Math. 22 (2004), 385–393. \url{https://doi.org/10.1016/S0723-0869(04)80016-1}
\bibitem{engelking} 
R. Engelking. \textit{General Topology}. Heldermann, Berlin, 1989.
\bibitem{Hu}
S. T. Hu. \textit{Theory of Retracts}. Wayne State University Press, Detroit, 1965.
%\bibitem{hfp}
%I. M. James. \textit{General Topology and Homotopy Theory}. Springer-Verlag, New York, 1984. 
%\bibitem{j1}
%J. Jaworowski. \textit{Equivariant extension of maps}. Pac. J. Math. 45 (1973), 229-244.  \url{https://doi.org/10.2140/pjm.1973.45.229}
%\bibitem{j2}
%J. Jaworowski. \textit{Extension of G-maps and Euclidean G-retracts}. Math. Z. 146 (1976), 143-148. \url{https://doi.org/10.1007/BF01187702}
%\bibitem{j3}
%J. Jaworowski. \textit{Extension properties of G-maps}. Proc. Intern. Conf. Geom. Topology, 209-213. PWN-Polish Sci. Publ., Warsaw, 1980.
\bibitem{hugo}
H. Juárez-Anguiano. \textit{Equivariant retracts and a topological social choice model}. Topol. Appl. 279 (2020). \url{https://doi.org/10.1016/j.topol.2020.107246}
\bibitem{hugo2} H. Juárez-Anguiano. \textit{Social choice models and topological dynamics}. Bol. Soc. Mat. Mex.  30(2) (2024) 62.  \url{https://doi.org/10.1007/s40590-024-00640-5}

%\bibitem{k} A. Kolmogorov. \textit{Sur la notion de moyenne}. Rendiconti Acad. Lincei 12 (1930), 388-391.
%\bibitem{sakai}
%K. Sakai. \textit{Geometric Aspects of General Topology}. Springer, Tokyo, 2013.
%\bibitem{smirnov}
%J. M. Smirnov. \textit{Sets of H-fixed points are absolute extensors}. Math. USSR-Sb. 27 (1975), 85-92.  \url{https://doi.org/10.1070/SM1975v027n01ABEH002501}
\bibitem{vanmill}
J. van Mill. \textit{Infinite-Dimensional Topology. Prerequisites and Introduction}. North-Holland, Amsterdam, 1989.
\bibitem{hip}
S. Weinberger. \textit{On the topological social choice model}. J. Econ. Theory 115 (2004), 377-384. \url{https://doi.org/10.1016/S0022-0531(03)00257-6}
\bibitem{we}
J. West. \textit{Induced involutions on Hilbert cube hyperspaces}. In Topology Proceedings, Vol. I (Conf., Auburn Univ., Auburn, Ala., 1976) (1977), 281-293.
\bibitem{ww} J. West and R. Wong. \textit{Based-free actions of finite groups on Hilbert cubes with absolute retract orbit spaces are conjugate}. In Geometric Topology (Proc. Georgia Topology Conf., Athens, Ga., 1977) (1979), Academic Press, New York - London, 655-672. \url{https://doi.org/10.1016/B978-0-12-158860-1.50043-5}

\bibitem{willard} 
S. Willard. \textit{General Topology}. Dover Publications, Inc., New York, 2004.
\bibitem{wo}
R. Wong. \textit{Periodic actions on the Hilbert cube}. Fund. Math. 85 (1974), 203-210.  \url{https://doi.org/10.4064/FM-85-3-203-210}
\end{thebibliography}
\end{document}